\documentclass[centertags,12pt]{amsart}
\usepackage{latexsym}
\usepackage{amssymb}

\numberwithin{equation}{section}
\makeatletter
\renewcommand{\subsection}{\@startsection
{subsection}{2}{0mm}{\baselineskip}{-0.25cm}
{\normalfont\normalsize\bf}}
\makeatother

\newtheorem{theorem}{Theorem}[section]

\newtheorem{lemma}[theorem]{Lemma}
\newtheorem{claim}[theorem]{Claim}
{\theoremstyle{remark}
\newtheorem{remark}[theorem]{Remark}}
\theoremstyle{definition}

\newtheorem{example}[theorem]{Example}

\def\P{\mathbf P}
\def\N{\mathbf N}

\def\cK{\mathcal K}
\def\cX{\mathcal X}

\def\sq{\sqrt{q}}
\def\fq{\mathbf F_q}

\def\deg{{\rm deg}}
\def\det{{\rm det}}
\def\frx{{\mathbf Fr}_{\mathcal X}}

\begin{document}

\author[Giulietti]{M. Giulietti}
\author[Pambianco]{F. Pambianco}
\author[Torres]{F. Torres}
\author[Ughi]{E. Ughi}\thanks{1991 Math. Subj. Class.: 
Primary 05B, Secondary 14H}

\title[Arcs and curves]{On complete arcs arising from plane curves}
\address{Dipartimento di Matematica - Universit\`a degli Studi di Perugia,
Via Vanvitelli 1 - 06123 Perugia (Italy)}
\email{giuliet@dipmat.unipg.it}
\email{fernanda@dipmat.unipg.it}
\address{IMECC-UNICAMP, Cx. P. 6065, Campinas, 13083-970-SP, Brazil}
\email{ftorres@ime.unicamp.br}
\address{Dipartimento di Matematica - Universit\`a degli Studi di Perugia,
Via Vanvitelli 1 - 06123 Perugia (Italy)}
\email{ughi@dipmat.unipg.it}

      \begin{abstract} 
We show that the set of $\fq$-rational points of either certain Fermat
curves or certain $\fq$-Frobenius non-classical 
plane curves is a complete $(k,d)$-arc in $\P^2(\fq)$, where $k$ and $d$
are respectively the number of $\fq$-rational points and the degree of the 
underlying curve.
      \end{abstract}

\maketitle

\section{Introduction and statement of results}\label{s1} 

A {\em $(k,d)$-arc} $\cK$ in the projective plane $\P^2(\fq)$, $\fq$ being
the finite field with $q$ elements, is a set of $k$ elements such that no
line in $\P^2(\fq)$ meets $\cK$ in more than $d$ points. The $(k,d)$-arc
is called {\em complete} if it is not contained in a $(k+1,d)$-arc. For
basic facts on arcs the reader is refered to \cite[Ch. 12]{h} (see also
the references therein), \cite[Sec. 5]{hirschfeld-storme},  
\cite{ball-blokhuis}, and \cite{szonyi}. 

A natural example of a $(k,d)$-arc is the set $\cX(\fq)$ of $\fq$-rational
points of a plane curve $\cX$ without linear components and defined over
$\fq$, where $k=\#\cX(\fq)$ and $d$ is the degree of $\cX$. As a matter of
terminology, we shall say that $\cX$ has the {\em arc property} whenever
$\cX(\fq)$ is a complete $(k,d)$-arc with $k$ and $d$ as above. As a
matter of fact, the interplay between the theory of algebraic curves 
and finite geometries was initiated by Segre around 1955. In \cite{segre}
(see also \cite[Sec. 10.4]{h}) he established an upper bound for the
second largest size that a complete $(k,2)$-arc in $\P^2(\fq)$ can have.
He proved this result by applying the Hasse-Weil upper bound to the
non-singular model of the envelope associated to $(k,2)$-arcs. For further
results on these arcs see Hirschfeld and Korchm\'aros' papers \cite{hk1}
and \cite{hk2}. For other applications of curves to finite geometries see
the surveys \cite{sz} and \cite{sz-sz}.

In this paper we are concerning with the problem of determining plane
curves having the arc property. This was asked around 1988 by Hirschfeld
and Voloch \cite[Problem III]{hirschfeld-voloch}. Only few examples of
such curves are known. Among them we have the irreducibles conics in odd
characteristic \cite[Ch. 8]{h}, certain cubics \cite[Sect.
5]{hirschfeld-voloch}, \cite[Sect. 6]{giulietti}, and Hermitian curves
\cite[Lemma 7.20]{h}. 

For $q$ a square, the Hermitian curve $\cX_{\sq+1}$ can be defined by
  \begin{equation}\label{eq1.1}
X^{\sq+1}+Y^{\sq+1}+Z^{\sq+1} = 0\, .
  \end{equation}
Then $\cX_{\sq+1}(\fq)$ is a complete $(q\sq+1,\sq+1)$-arc (loc. cit.). 
The completeness property means that for each $P\in\P^2(\fq)\setminus
\cX_{\sq+1}(\fq)$ 
there exists a $\fq$-rational line $\ell$ such that $\# \ell\cap
\cX_{\sq+1}(\fq) = \sq+1$. This property can be easily shown by using the
special feature of Eq. (\ref{eq1.1}) (see the proof of Theorem
\ref{thmA}(1)). Moreover, this property is also a 
consequence of the fact that the image of $P\in \cX_{\sq+1}$ by the
$\fq$-Frobenius morphism lies on the tangent line of $\cX_{\sq+1}$ at $P$ 
(see the proof of Theorem \ref{thmB}). A plane 
curve satisfying the above property for general points is called 
{\em $\fq$-Frobenius non-classical} \cite{sv}, \cite{hefez-voloch}.

The Hermitian curve is a member of the family of Fermat curves 
$\cX_d(a,b)$ defined by 
   \begin{equation}\label{eq1.2}
aX^d+bY^d+Z^d=0\, ,
   \end{equation}
where $d=(q-1)/(q'-1)$, $q=p^n$, $q'=p^e$, $p:={\rm
char}(\fq)$ such that $e<n$ and $e|n$, and where $a,b\in \mathbf
F_{q'}^*:=\mathbf F_{q'}\setminus\{0\}$.

The first aim of this paper is to extend the arc property of 
the Hermitian curve to the curves $\cX_n(a,b)$ above as well as to the
Fermat curve $\cX_{q-1}$ defined by
   \begin{equation}\label{eq1.3}
X^{q-1}+Y^{q-1}=2Z^{q-1}\, ,
   \end{equation}
provided that $p\ge 3$.
   \begin{theorem}\label{thmA} For $p\ge 3$, the following statements
hold:
  \begin{enumerate}
\item[\rm(1)] the set of $\fq$-rational points of the
Fermat curve defined by Eq. (\ref{eq1.2}) is a complete
$(k,d)$-arc, where $k=d(q-d+2);$
\item[\rm(2)] the set of $\fq$-rational points of the Fermat curve defined
by Eq. (\ref{eq1.3}) is a complete $((q-1)^2,q-1)$-arc.
\end{enumerate}
  \end{theorem}
The arc in part (2) of this theorem is maximal among 
$(k,q-1)$-arcs for which there exists an external line, see Remark
\ref{rem3.1}. On the other hand, it seems that the arcs in part (1) are
new.

We notice that the curves in Theorem \ref{thmA} are among the Fermat
curves having a large number of $\fq$-rational points
\cite{garcia-voloch1}. We also notice that the hypothesis $p\ge 3$ is
necessary, see Remark \ref{rem3.0}. 

The second aim of this paper is to show that certain 
$\fq$-Frobenius non-classical plane curves do satisfy the arc property.
   \begin{theorem}\label{thmB} Let $\cX$ be a non-singular
$\fq$-Frobenius non-classical plane curve of degree $d$. Let $\epsilon$ be
the order of contact of $\cX$ with the tangent at a general point. If 
  \begin{equation}\label{eq1.4}
d(d-1)<(q+1)\epsilon\, ,
  \end{equation}
then the set of $\fq$-rational points of $\cX$ is a complete
$(k,d)$-arc, where $k=\#\cX(\fq)=d(q-d+2)$.
  \end{theorem}
The Fermat curve in Theorem \ref{thmA}(2) is $\fq$-classical by 
\cite[Thm. 2]{garcia-voloch1}. Hence the hypothesis of being
$\fq$-Frobenius non-classical in Theorem \ref{thmB} is not necessary. We
observe that $d(d-1)/\epsilon$ is the degree of the dual curve of $\cX$;
see Remark \ref{rem2.3}.

For the Hermitian curve $\cX_{\sq+1}$, $\epsilon=\sq$ (see e.g.
\cite{garcia-viana}) and therefore it satisfies (\ref{eq1.4}) in
Theorem \ref{thmB}. The Fermat curves in Theorem \ref{thmA}(1) are
$\fq$-Frobenius non-classical, see \cite[Thm.
2]{garcia-voloch1}. For these curves, $\epsilon=q'$ \cite[Thm.
2]{hefez-voloch}, and so they satisfy 
(\ref{eq1.4}) if and only if they are Hermitian curves. So far, we
could not find examples of non-singular $\fq$-Frobenius non-classical curves 
fulfilling (\ref{eq1.4}) which are not $\fq$-isomorphic to Hermitian
curves; see Remark \ref{rem4.1}. 

\section{Frobenius non-classical planes curves}\label{s2}

The study of Frobenius non-classical curves was initiated by Hefez and
Voloch \cite{hefez-voloch} based on a fundamental paper by St\"ohr and
Voloch \cite{sv}, where an approach to the Hasse-Weil bound was given.

In this paper we only consider irreducible non-linear plane   
curves defined over $\fq$. Let $\cX\subseteq \P^2(\bar\fq)$ be such a
curve. For $i=0,1,2$, let $x_i$ be the coordinates functions 
of $\P^2(\bar\fq)$ on $\cX$. Let $t$ be a separating variable
of $\fq(\cX)|\fq$ and denote by $D^i=D^i_t$ the $i$-th Hasse
derivative on $\cX$. The order sequence of $\cX$ (see \cite[p. 5]{sv}) are
the numbers $0,1$ and $\epsilon=\epsilon(\cX)$, where $\epsilon>1$ is the
least integer such that
$$
\det\begin{pmatrix} x_0 & x_1 & x_2\\
              D^1x_0 & D^1x_1 & D^1x_2\\
              D^\epsilon x_0 & D^\epsilon x_1 & D^\epsilon x_2
\end{pmatrix}\neq 0\, .
$$
Geometrically, the numbers 0, 1 and $\epsilon$ represent all the possible
intersection multiplicities of the curve $\cX$ with lines in 
$\P^2(\bar\fq)$ at general points.

The $\fq$-Frobenius order sequence of $\cX$ (see \cite[p. 9]{sv}) are the
numbers 0 and $\nu=\nu(\cX,q)$, where $\nu>0$ is the least integer such
that 
  \begin{equation}\label{eq2.1}
\det\begin{pmatrix} x_0^q & x_1^q & x_2^q\\
               x_0   & x_1   & x_2\\
               D^\nu x_0 & D^\nu x_1 & D^\nu x_2
  \end{pmatrix}\neq 0\, .
   \end{equation}
We have that $\nu\in \{1,\epsilon\}$ \cite[Prop. 2.1]{sv}. The plane 
curve 
$\cX$ is called $\fq$-Frobenius non-classical if $\nu=\epsilon$ (or
equivalently if $\nu>1$). 
   \begin{remark}\label{rem2.1} From Eq. (\ref{eq2.1}) follows that 
$\cX$ is $\fq$-Frobenius non-classical if and only if $\frx(P)\in T_P\cX$
for all non-singular points $P\in \cX$, where $\frx$ is the
$\fq$-Frobenius 
morphism on $\cX$ and $T_P\cX$ is the tangent line to $\cX$ at $P$.
   \end{remark}
   \begin{remark}\label{rem2.2} Let $\cX$ be a curve as above.

(i) If $\epsilon(\cX)>2$ then it is a power of $p$
\cite[Prop. 2]{garcia-voloch}.

(ii) If $\cX$ is $\fq$-Frobenius non-classical and $p>2$, then
$\epsilon(\cX)>2$ \cite[Prop. 1]{hefez-voloch}.

(iii) If $\cX$ is $\fq$-Frobenius non-classical, then $\epsilon(\cX)\le q$
\cite[p. 266]{hefez-voloch}. If in addition $\cX$ is non-singular, then
$\epsilon(\cX)\le \sqrt{q}$ \cite[Prop. 6]{hefez-voloch}.

(iv) Let $\cX$ be non-singular $\fq$-Frobenius non-classical plane curve. 
Let $d$ the degree of $d$ and suppose that $\epsilon=\epsilon(\cX)>2$. 
Then \cite[Props. 5, 6]{hefez-voloch}
$$
\sqrt{q}+1\le d\le (q-1)/(\epsilon-1)\, .
$$
    \end{remark}
    \begin{remark}\label{rem2.3} Let $\cX$ be a non-singular plane curve
of degree $d\ge 2$, 
$\cX^*$ the dual curve of $\cX$ and $T:\cX\to \cX^*$ the dual map; i.e,
$T(P)=T_P\cX$. Then $d(d-1)=\deg(T)d^*$, where $d^*$ is the degree of
$\cX^*$; this follows e.g. from \cite[Lemma 4.3]{homma}. Now by a result
of Kaji \cite[Cor. 4.5]{kaji}, $T$ is purely inseparable. If in addition,
$\cX$ is $\fq$-Frobenius non-classical, then $\deg(T)=\epsilon(\cX)$
\cite[Props. 3 and 4]{hefez-voloch}. Therefore, in Theorem \ref{thmB} we
look for $\fq$-Frobenius non-classical plane curves such that the degree
of its dual curve is upper bounded by $(q+1)$.
    \end{remark}

Let $F=F(X_0,X_1,X_2)=0$ be the equation of $\cX$ over $\fq$. From 
\cite[Thm. 1]{garcia-voloch},  
$\epsilon=\epsilon(\cX)>2$ if and only if there exist homogeneous
polynomials $H,
P_0, P_1, P_2\in \fq[X_0,X_1,X_2]\setminus\{0\}$ such that
  \begin{equation}\label{eq2.2}
FH=X_0P_0^\epsilon+X_1P_1^\epsilon+X_2P_2^\epsilon\, .
  \end{equation}
We have that $H=1$ if $\cX$ is non-singular (loc. cit, p. 462). Now 
it is easy to see that $\cX$ is $\fq$-non-classical if and only if there 
exists $H_1\in\fq[X_0,X_1,X_2]$ such that
  \begin{equation}\label{eq2.3}
FH_1=X_0^{q/\epsilon}P_0+X_1^{q/\epsilon}P_1+X_2^{q/\epsilon}P_2\, .
  \end{equation}
Finally, we mention a formula for the precise number of $\fq$-rational
points of non-singular $\fq$-Frobenius non-classical curves.
   \begin{lemma}\label{lemma2.1} {\rm (\cite[Thm. 1]{hefez-voloch})} 
Let $\cX$ be a plane non-singular $\fq$-Frobenius non-classical curve of
degree $d$. Then 
$$
\#\cX(\fq)=d(q-d+2)\, .
$$
   \end{lemma}

\section{Proof of Theorem \ref{thmA}}\label{s3}

(1) That $\#\cX_d(a,b)=d(q-d+2)$ is
well known, see e.g. \cite[p. 354]{garcia-voloch1}. (This 
result also follows from \cite[Thm. 2]{garcia-voloch1} and Lemma
\ref{lemma2.1}). 

Next, for $P\in \P^2(\fq)\setminus \cX_d(\fq)$, we will show that there
exists a 
$\fq$-rational line $\ell$ which passes through $P$ and intersects the
curve $\cX_d(a,b)$ in $d$ distinct $\fq$-rational points. We recall the
following easy fact.
  \begin{claim}\label{claim3.1} Let $A,B\in \mathbf F_{q'}^*$. Then the
equation $AX^d+B=0$ has $d$ distinct solutions in $\fq$.
   \end{claim}
  \begin{proof} Since $p$ does not divide $d$, the equation has $d$
solutions in 
$\bar\fq$. If $x$ is a solution, then $x^{q-1}=1$, as $d=(q-1)/(q'-1)$,
and hence $x\in \fq$.
  \end{proof}
We consider four cases:

{\bf Case 1:} $P=(\alpha:\beta:0)$. Let $\ell: Z=0$. Then $P\in \ell$ and
$$ 
\ell\cap\cX_d(a,b)=\{(\lambda:1:0): a\lambda^d+b=0\}\, .
$$
{\bf Case 2:} $P=(\alpha:\beta:1)$ and $a\alpha^d+1\neq 0$. Let $\ell: 
X=\alpha Z$. Then $P\in \ell$ and
$$
\ell\cap\cX_d(a,b)=\{(\alpha:\lambda:1): b\lambda^d+a\alpha^d+1=0\}\, .
$$
{\bf Case 3:} $P=(\alpha:\beta:1)$ and $b\beta^d+1\neq 0$. Let $\ell: 
Y=\beta Z$. Then $P\in \ell$ and
$$
\ell\cap\cX_d(a,b)=\{(\lambda:\beta:1): a\lambda^d+b\beta^d+1=0\}\, .
$$
{\bf Case 4:} $P=(\alpha:\beta:1)$ and $a\alpha^d+1=b\beta^d+1=0$. Let
$\ell: \beta X=\alpha Y$. Then $P\in \ell$ and
$$
\ell\cap\cX_d(a,b)=\{(\alpha,\beta,\beta\lambda): \lambda^d+2b=0\}\, .
$$
Now by Claim \ref{claim3.1} all the sets above have cardinality $d$ and  
are contained in $\cX_d(a,b)(\fq)$. This completes the proof of Theorem
\ref{thmA}(1).

(2) We have that
$$
\cX_{q-1}(\fq)=\{(\alpha:\beta:\gamma)\in \P^2(\fq): \alpha\beta\gamma\neq
0\}\, ,
$$
so that $\#\cX_{q-1}(\fq)=q^2+q+1-3q=(q-1)^2$. Let 
$P=(\alpha:\beta:\gamma)\in \P^2(\fq)\setminus \cX_{q-1}(\fq)$.
The proof of the arc property for $\cX_{q-1}$ follows from the following
five computations.

{\bf Case 1:} $\beta=\gamma=0$. Let $\ell: Y=Z$. Then $P\in \ell$ and
$$
\ell\cap\cX_{q-1}=\{(\lambda:1:1): \lambda\in \fq^*\}\, .
$$
{\bf Case 2:} $\beta\neq 0$ and $\gamma=0$. Let $\ell: Y=bX$. Then $P\in
\ell$ and
$$
\ell\cap\cX_{q-1}=\{(1:b:\lambda): \lambda\in \fq^*\}\, .
$$
{\bf Case 3:} $\alpha=\beta=0$. Let $\ell: X=Y$. Then
$P\in \ell$ and
$$
\ell\cap\cX_{q-1}=\{(1:1:\lambda): \lambda\in \fq^*\}\, .
$$
{\bf Case 4:} $\alpha\neq 0$ and $\gamma=1$. Let $\ell: X=aZ$. Then $P\in
\ell$ and
$$
\ell\cap\cX_{q-1}=\{(a:\lambda:1): \lambda\in \fq^*\}\, .
$$
{\bf Case 5:} $\beta\neq 0$ and $\gamma=1$. Let $\ell: Y=bZ$. Then $P\in
\ell$ and 
$$
\ell\cap\cX_{q-1}=\{(\lambda:b:1): \lambda\in \fq^*\}\, .
$$
   \begin{remark}\label{rem3.0} The hypothesis $p\ge 3$ in Theorem
\ref{thmA}(1) is necessary. Indeed, consider the Fermat curve $\cX$ 
defined by 
$$
X^{q-1}+Y^{q-1}+Z^{q-1}=0\, ,
$$
where $q$ is a power of two. Take $P=(1:1:1)$. Then it is easy to see that
$\#\ell\cap\cX<q-1$ for any $\fq$-rational line $\ell$ passing through
$P$.
   \end{remark}
   \begin{remark}\label{rem3.1} For a $(k,q-1)$-arc $\cK$ having 
an external line (i.e. a $\fq$-rational line $\ell$
such that $\ell\cap \cK=\emptyset$), we have that $k\le (q-2)q+1=(q-1)^2\
(*)$ \cite[Thm. 12.40]{h}. Since the line $\ell: Z=0$ is external to the  
$((q-1)^2,q-1)$-arc in Theorem \ref{thmA}(2), we have that 
the upper bound in $(*)$ is attained by this arc.
   \end{remark}
   \begin{remark}\label{rem3.2} Notice that the $((q-1)^2,q-1)$-arc in
Theorem \ref{thmA}(2) is the complement of three non-concurrent lines in
$\P^2(\fq)$. Moreover, it can be shown that any such set is a complete
$((q-1)^2,q-1)$-arc which is also the set of $\fq$-rational points of a 
plane curve of degree $q-1$; it turns out that this curve is
$\fq$-isomorphic to the Fermat curve $\cX_{q-1}$.
    \end{remark}
   \begin{example}\label{ex3.1} Taking $q=p^3$ and $q'=p$ in Theorem
\ref{thmA}(1), we have that there exists a complete $(k,p^2+p+1)$-arc in
$\P^2(\fq)$ with $k=(p^2+p+1)(p^3-p^2-p+1)$. In particular, there exists a
complete $(208, 13)$-arc in $\mathbf F_{27}$.
   \end{example}

\section{Proof of Theorem \ref{thmB}}\label{s4}

Theorem \ref{thmB} will be a consequence of the following more general
result.
   \begin{theorem}\label{thm4.1} Let $\cX$ be a $\fq$-Frobenius
non-classical (possible singular) plane curve. Let $d$ be the degree of
$\cX$ and $k=\#\cX(\fq)$. Then $\cX(\fq)$ is a complete $(k,d)$-arc
provided that
$$
k>(d-\epsilon)(q+1-\#S)+(d-1)\#S \, ,
$$
where $S$ is the set of singular points of $\cX$ and $\epsilon$ is as in 
Theorem \ref{thmB}.
   \end{theorem}
   \begin{proof} Suppose that $\cX(\fq)$ is not complete. Then there
exists $P\in 
\P^2(\fq)\setminus \cX(\fq)$ such that for any $\fq$-rational line $\ell$
through $P$,
   \begin{equation}\label{eq4.1}
\#\ell\cap\cX(\fq)<d\, .
   \end{equation}
   \begin{claim}\label{claim4.1} If $\ell$ is a line such that
(\ref{eq4.1}) holds and that $\ell\cap S=\emptyset$, then
$$
\#\ell\cap\cX(\fq)\le (d-\epsilon)\, .
$$
  \end{claim}
  \begin{proof} {\bf Case 1: $\ell\cap \cX\not\subseteq \cX(\fq)$.} Let
$Q\in 
\ell\cap\cX\setminus \cX(\fq)$. Then, as $\ell$ is $\fq$-rational,
$\frx(Q)\in \ell$ and thus, as $\ell\cap S=\emptyset$, $\ell$ is the
tangent line of $\cX$ at $Q$ (see 
Remark \ref{rem2.1}). Therefore
$$
\#\ell\cap\cX(\fq)\le d-\epsilon\, ,
$$
since $I(\cX,\ell;Q)\ge \epsilon$; see \cite[p.5]{sv}.

{\bf Case 2: $\ell\cap\cX\subseteq \cX(\fq)$.} From (\ref{eq4.1}) and 
Bezout's theorem there exists $Q\in \ell\cap\cX(\fq)$ such that
$j(Q):=I(\cX,\ell;Q)>1$. Then, as $\ell\cap S=\emptyset$, $\ell$ is the
tangent line of $\cX$ at $Q$. We have that 
$$
\#\ell\cap\cX(\fq)\le d-j(Q)+1
$$
and the claim follows from the fact that $j(Q)\ge \epsilon+1$ \cite[Cor.
2.10]{sv}.
   \end{proof}
Now there are at most $N\le \# S$ lines $\ell^\prime$ such that
$\ell^\prime \cap S\neq \emptyset$. For each of these lines, 
$\ell^\prime\cap\cX(\fq)\le d-1$; hence from Claim \ref{claim4.1} we have
$$
k\le (d-\epsilon)(q+1-N)+(d-1)N\, .
$$
Then $k\le (d-\epsilon)(q+1)+N(\epsilon-1)\le
(d-\epsilon)(q+1)+(\epsilon-1)\#S$, a contradiction. This finishs the  
proof of Theorem \ref{thm4.1}.
    \end{proof}
{\em Proof of Theorem \ref{thmB}.} We have $S=\emptyset$ and so the
hypothesis in Theorem \ref{thm4.1} is 
\begin{equation}\label{eq4.2}
k>(d-\epsilon)(q+1)\, .
\end{equation}
Since $k=d(q-d+2)$ (see Lemma \ref{lemma2.1}), it turns out that
(\ref{eq4.2}) is equivalent to $d(d-1)<(q+1)\epsilon)$ and the result
follows.
    \begin{remark}\label{rem4.1} Let $\cX$ be a
non-singular $\fq$-Frobenius non-classical curve of degree $d$. Assume 
$\epsilon=\epsilon(\cX)>2$. Then from Eqs. (\ref{eq2.2}) and
(\ref{eq2.3}), $d=\lambda\epsilon+1$ for some $\lambda \in \N$. If 
(\ref{eq1.4}) holds, then
$$
\sq/\epsilon \le \lambda < \sqrt{q/\epsilon}\, ,
$$
where the first inequality follows from Remark \ref{rem2.2}(iv). 

For a concrete example take $q=p^3$. Then $\epsilon=p$ by Remark
\ref{rem2.2}(i)(iii) and so $\sqrt{p}\le \lambda<p$. Therefore $\cX$ will
satisfy (\ref{eq1.4}) if $\lambda=p-1$, i.e. if $\cX$ has
degree $d=p^2-p+1$. The existence of a such curve is equivalent to the
existence of polynomials $P_0, P_1, P_2\in \fq[X_0,X_1,X_2]$ of degree
$p-1$ and $H_1\in \fq[X_0,X_1,X_2]$ such that Eqs.
(\ref{eq2.2}) and (\ref{eq2.3}) hold true. This seems
an involved problem.  On the other hand, the curve $\cX$ will rise to the
existence of a complete $((p^2-p+1)(p^3-p^2+p+1), p^2-p+1)$-arc.
Unfortunately, the existence of such an arc is not known. 
    \end{remark}

{\bf Acknowledgments.} The authors wish to thank J.W.P. Hirschfeld and G.
Korchm\'aros for useful comments. This research was carried out with the
support of the Italian Ministry for Research and Technology (project 40\%
``Strutture geometriche, combinatorie e loro applicazioni"). Part of this
paper was written while Torres was visiting Perugia in February and
December 1999.


\begin{thebibliography}{99}

\bibitem{ball-blokhuis} S. Ball and A. Blokhuis, On the incompleteness of
$(k,n)$-arcs in Desarguesian planes of order $q$ where $n$ divides $q$,
{\em Geom. Dedicata} {\bf 74} (1999), 325--332.

\bibitem{garcia-viana} A. Garcia and P. Viana, Weierstrass points on
certain non-classical curves, {\em Arch. Math.} {\bf 46} (1986), 315--322.

\bibitem{garcia-voloch} A. Garcia and J.F. Voloch, Wronskians and
linear independence in fields of prime characteristic, {\em Manuscripta
Math.} {\bf 59} (1987), 457--469.

\bibitem{garcia-voloch1} A. Garcia and J.F. Voloch, Fermat curves over
finite fields, {\em J. Number Theory} {\bf 30} (1988), 345--356.

\bibitem{giulietti} M. Giulietti, On plane arcs and cubic curves,
preprint.

\bibitem{hefez-voloch} A. Hefez and J.F. Voloch, Frobenius non
classical curves, {\em Arch. Math.} {\bf 54} (1990), 263--273.

\bibitem{h} J.W.P. Hirschfeld, {\em Projective Geometries Over Finite
Fields}, second edition, Oxford University Press, Oxford, 1998.

\bibitem{hk1} J.W.P. Hirschfeld and G. Korchm\'aros, On the number of
rational points on an algebraic curve over a finite field, {\em Bull.
Belg. Math. Soc. Simon Stevin} {\bf 5} (1998), 313--340.

\bibitem{hk2} J.W.P. Hirshfeld and G. Korchm\'aros, Arcs and curves over a
finite field, {\em Finite Fields Appl.} {\bf 5} (1999), 393--408.

\bibitem{hirschfeld-storme} J.W.P. Hirschfeld and L. Storme, The
packing problem in statistics, coding theory and finite projective spaces,
{\em J. Statist. Plann. Inference} {\bf 72} (1998), 355--380.

\bibitem{hirschfeld-voloch} J.W.P. Hirschfeld and J.F. Voloch, The
characterization of elliptic curves over finite fields, {\em J. Austral. 
Math. Soc. Ser. A} {\bf 45} (1988), 275--286.

\bibitem{homma} M. Homma, Funny curves in characteristic $p>0$, {\em Comm.
Algebra} {\bf 15}(7) (1987), 1469--1501.

\bibitem{kaji} H. Kaji, On the Gauss maps of space curves in
characteristic $p$, {\em Compositio Math.} {\bf 70} (1989), 177--197.

\bibitem{segre} B. Segre, Ovals in a finite projective plane, {\em Canad.
J. Math.} {\bf 7} (1955), 414--416.

\bibitem{sv} K.O. St\"ohr and J.F. Voloch, Weierstrass points and
curves over finite fields, {\em Proc. London Math. Soc.} {\bf 52} (1986),
1--19.

\bibitem{sz-sz} P. Sziklai and T. Sz\"onyai, Blocking sets and algebraic
curves, {\em Rend. Circ. Mat. Palermo Suppl.} {\bf 51} (1998), 71--86.

\bibitem{sz} T. Sz\"onyi, Some applications of algebraic curves in finite
geometry and combinatorics, ``Surveys in Combinatorics" (R.A. Bailey Ed.),
197--236, Cambridge Univ. Press, Cambridge, 1997.

\bibitem{szonyi} T. Sz\"onyi, On the embedding of $(k,p)$-arcs in maximal
arcs, {\em Design, Codes and Cryptography} {\bf 18} (1999), 235--246.

\end{thebibliography}
\end{document}